\documentclass[11pt]{article}
\usepackage{amsfonts,amsmath}
\usepackage{latexsym}
\usepackage{amsthm}
\usepackage{amscd}
\usepackage{commath}
\usepackage{epsfig}
\usepackage{amssymb}
\usepackage{caption}
\usepackage{enumitem}
\usepackage{mathtools}
\usepackage{booktabs}
\usepackage[table]{xcolor}
\usepackage[toc,page]{appendix}
\usepackage[english]{babel}

\usepackage[letterpaper,top=2cm,bottom=2cm,left=2.5cm,right=2.5cm,marginparwidth=1.5cm]{geometry}

\usepackage[colorlinks=true, allcolors=blue,bookmarksnumbered]{hyperref}
\usepackage{color,soul}
\usepackage{float,subcaption}
\usepackage{graphicx}
\usepackage{stackengine}
\usepackage{tikz}
\usetikzlibrary{intersections}
\usepackage[export]{adjustbox}
\usepackage{array}
\usepackage[percent]{overpic}

\usepackage[capitalise]{cleveref}

\usepackage{lineno}

\newtheorem{defn}{Definition}[section]

\newtheorem{theorem}[defn]{Theorem}

\newtheorem{lem}[defn]{Lemma}
\newtheorem{cor}[defn]{Corollary}
\theoremstyle{remark}
\newtheorem{remark}[defn]{Remark}

\numberwithin{equation}{section}
\numberwithin{figure}{section}

\newcommand{\bb}{\begin{equation}}
\newcommand{\ee}{\end{equation}}

\newcommand{\MD}[1]{{\color{purple} [Matteo: #1]}}
\newcommand{\LC}[1]{{\color{blue}#1}}
\newcommand{\ALN}[1]{{\color{orange}#1}}

  \newcommand{\newALN}[1]{{#1}}

\newcommand{\origin}{\mathbf{o}}

\newcolumntype{?}{!{\vrule width 1.5pt }}
\newlength\savedwidth

\newcommand{\tightoverset}[2]{%
  \mathop{#2}\limits^{\vbox to -.5ex{\kern-0.75ex\hbox{$#1$}\vss}}}

\newcommand\sut{\,;\ }
\newcommand{\dfn}[1]{\textbf{\textit{#1}}}
\newcommand\HH{\mathbb{H}} 

\newcommand{\R}{\mathbb{R}}

\newcommand{\IC}{\mathrm{IC}_{p,q}} 

\def\rlabel #1 #2{\begin{equation} \label{#1} #2 \end{equation}}

\def\rproof{\begin{proof}}

\def\Qed{\end{proof}}

\title{{\Large
	 \textsc{Extremal Ising Gibbs States on Lobachevsky lattices}
	}}
\author{
Matteo \textsc{D'Achille}\thanks{Laboratoire de Math\'ematiques d'Orsay, CNRS, Universit\'e Paris-Saclay, F-91405, Orsay, France\newline $_{}$\hfill    \href{mailto:matteo.dachille@universite-paris-saclay.fr}
{\texttt{matteo.dachille@universite-paris-saclay.fr}}},\hspace{2pt}
Loren \textsc{Coquille}\thanks{ Institut Fourier, Universit\'e Grenoble Alpes \& CNRS, F-38000 Grenoble, France\newline $_{}$\hfill    \href{mailto:loren.coquille@univ-grenoble-alpes.fr}
{\texttt{loren.coquille@univ-grenoble-alpes.fr}}},\hspace{2pt}
Arnaud \textsc{Le Ny}\thanks{Univ Paris Est Cr\'eteil \& CNRS,  LAMA UMR8050,  F-94010 Cr\'eteil, France \newline $_{}$\hfill     \href{mailto:arnaud.le-ny@u-pec.fr}
{\texttt{arnaud.le-ny@u-pec.fr}} } \hspace{5pt}
}
\date{\today}

\begin{document}
\maketitle

\begin{abstract}
{
	We exhibit an uncountable family of extremal {inhomogeneous} Gibbs measures of the low temperature Ising model on regular tilings of the hyperbolic plane. These states arise as low temperature perturbations of local ground states having a sparse enough set of frustrated edges, the sparseness being measured in terms of the isoperimetric constant of the graph.
	This result is implied by {an extension} of the article \cite{ckln} on regular trees to non-amenable graphs. We moreover argue how we can deduce  the extremality of an uncountable subset of the Series--Sinai states \cite{series1990ising} at low temperature.}
\medskip
\medskip

\end{abstract}
\newpage


\normalsize
\vspace{-0.5cm}

\section{Introduction}


Since a pioneering work by Lund--Rasetti--Regge~\cite{LRR77}, the ferromagnetic n.n.~Ising model on lattices {which are naturally} embedded in the hyperbolic plane $\HH_2$ has attracted considerable interest in the physics literature. {See also the work by  Rietman--Nienhuis--Oitmaa~\cite{RNO} or the more recent applied work in~\cite{BPR2020}.} 
We recall that on the Euclidean lattice $\mathbb{Z}^{2}$ (corresponding to the tiling $\mathcal{L}_{4,4}$ below), a celebrated result due to Aizenmann~\cite{aizenman1980} and Higuchi~\cite{higuchi1979} states that any Gibbs measure writes as a convex combination of the pure phases $\mu^+$ and $\mu^-$. On $\mathbb{Z}^{3}$, the famous result of Dobrushin \cite{Dob72} exhibits a countable family of extremal non translation-invariant states at low temperature, with a localised interface at a given height~\cite{9bis}.

When the underlying graph is a $d$-regular tree, the Ising model possesses uncountably many such interface states, extremal and non tree-automorphism invariant,  in the low temperature regime $\beta>\beta_c(d)$, as it was proved by Higuchi in the late seventies~\cite{hig77} (and in ~\cite{BLG91}). 
In 2008, Rozikov and Rakhmatullaev \cite{RR08} exhibit the so-called ``weakly periodic'' Gibbs measures,  which correspond to subgroups in the representation group of the Cayley tree. 
These states can be thought of as generalizations of Dobrushin states, but with many interfaces, possibly countably infinitely many, localised on a symmetric pattern.\\
{In 2012, Gandolfo, Ruiz and Shlosman \cite{GRS12,GRS15}, exhibit extremal inhomogeneous Gibbs measures arising
	as low temperature perturbations of {{ground states},  which have a sparse enough set of frustrated bonds}}. 
These states in general do not possess any symmetries of the tree. 
The two last authors, together with {K\"ulske}\ALN{,} rigorously proved in \cite{ckln} the statements of \cite{GRS12} and generalise them to a broader setting: non-compact state space gradient models, models without spin-symmetry, models in small random fields.

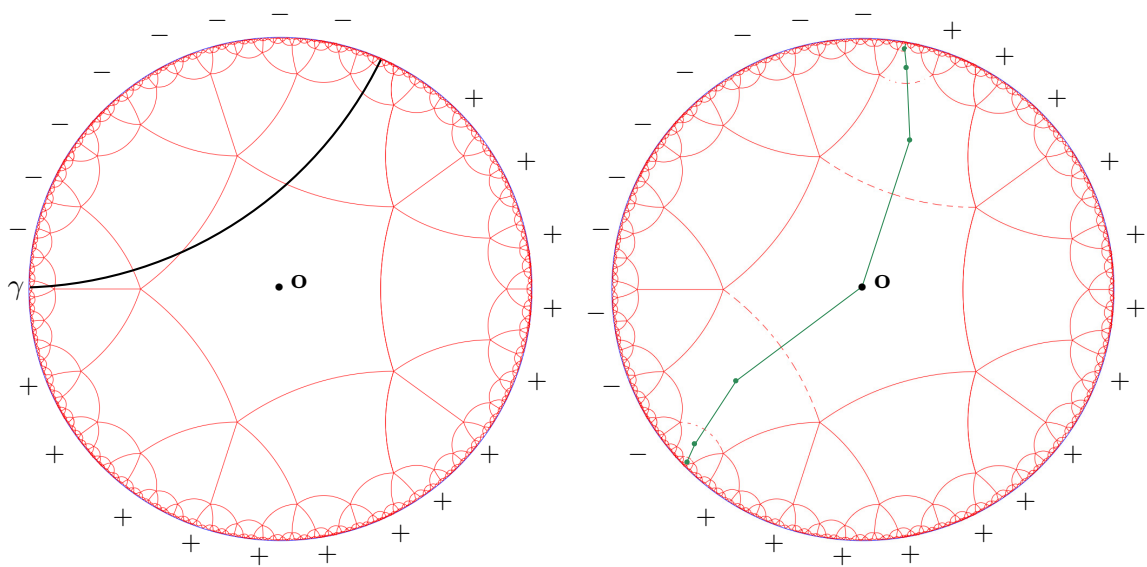
\begin{figure}[!htbp]
\centering
\begin{overpic}[width=.42\linewidth]{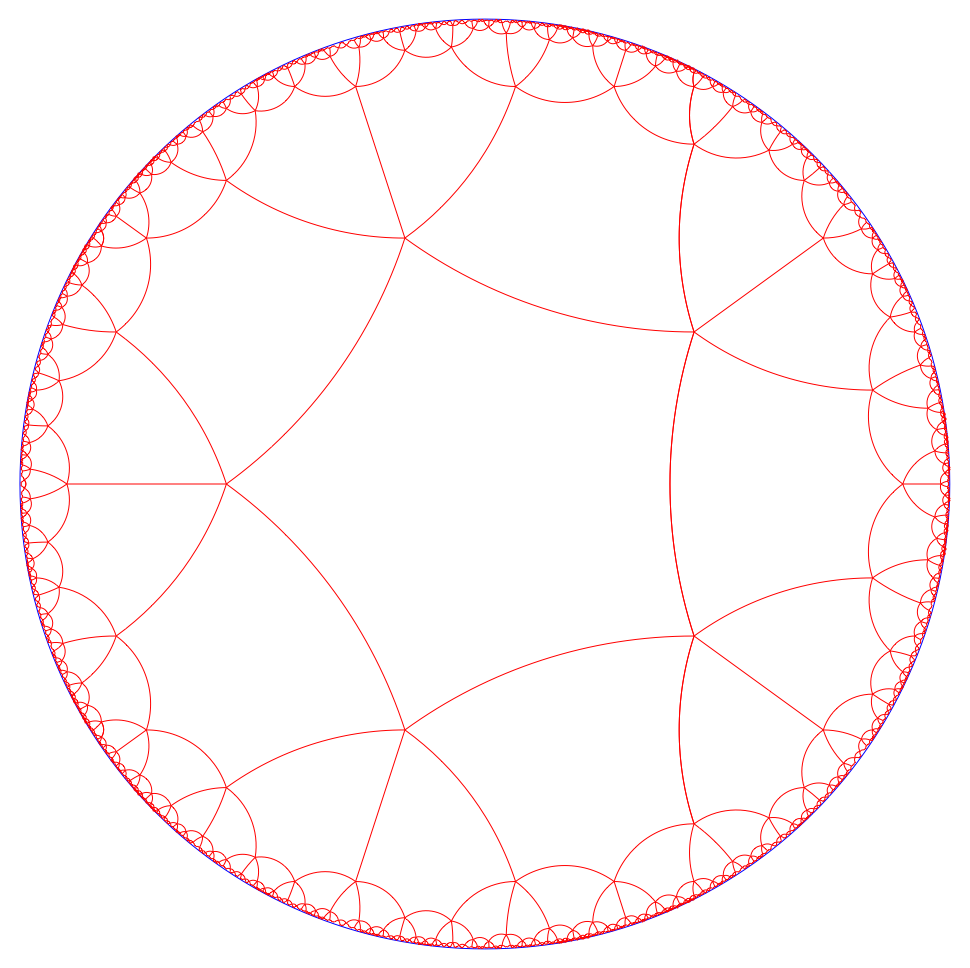}
\put (85,85) {{\small $+$}}
\put (95,73) {{\small $+$}}
\put (100,59) {{\small $+$}}
\put (100,45) {{\small $+$}}
\put (97,31) {{\small $+$}}
\put (88,17) {{\small $+$}}
\put (79,8) {{\small $+$}}
\put (71,2) {{\small $+$}}
\put (57,-2) {{\small $+$}}
\put (44,-2.5) {{\small $+$}}
\put (18,5) {{\small $+$}}
\put (5,17) {{\small $+$}}
\put (0,30) {{\small $+$}}
\put (31,0) {{\small $+$}}
\put (-2,60) {{\small $-$}}
\put (1,70) {{\small $-$}}
\put (6,80) {{\small $-$}}
\put (14,90) {{\small $-$}}
\put (25,97) {{\small $-$}}
\put (36,100) {{\small $-$}}
\put (48,101) {{\small $-$}}
\put (60,100) {{\small $-$}}
\put (52,50) {$\origin$}
\put (48.8,49.5) {{\tiny $\bullet$}}
\put(2,50){\tikz \draw[thick,black] (0,0) arc (-88:-26:5.4);}
\put (-2,49) {$\gamma$}
\end{overpic}
\hspace{15pt}
\begin{overpic}[width=.42\linewidth]{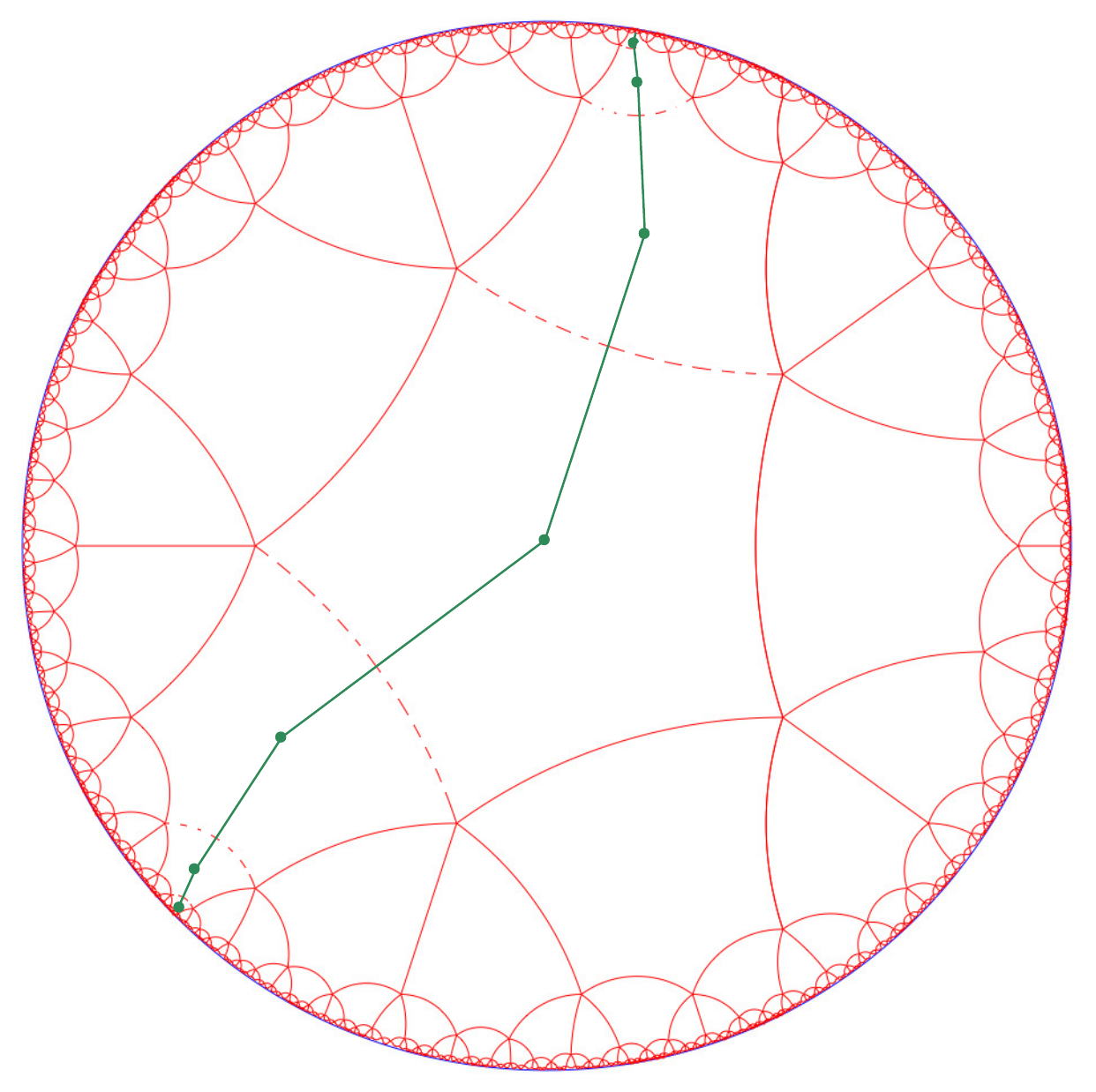}
\put (65,98) {{\small $+$}}
\put (75,93) {{\small $+$}}
\put (85,85) {{\small $+$}}
\put (95,73) {{\small $+$}}
\put (100,59) {{\small $+$}}
\put (100,45) {{\small $+$}}
\put (97,31) {{\small $+$}}
\put (88,17) {{\small $+$}}
\put (79,8) {{\small $+$}}
\put (71,2) {{\small $+$}}
\put (57,-2) {{\small $+$}}
\put (44,-2.5) {{\small $+$}}
\put (18,5) {{\small $+$}}
\put (5,17) {{\small $-$}}
\put (0,30) {{\small $-$}}
\put (31,0) {{\small $+$}}
\put (-2,60) {{\small $-$}}
\put (1,70) {{\small $-$}}
\put (6,80) {{\small $-$}}
\put (14,90) {{\small $-$}}
\put (25,97) {{\small $-$}}
\put (36,100) {{\small $-$}}
\put (48,101) {{\small $-$}}
\put (-3,44) {{\small $-$}}
\put (52,50) {$\origin$}
\put (48.8,49.5) {{\tiny $\bullet$}}
\end{overpic}

\caption{
A Series--Sinai state on $\mathcal{L}_{5,5}$ selected by $\gamma$, a bi-infinite geodesic of $\mathbb{H}_{2}$ (in black). Right: broken edges of $\mathcal{L}_{5,5}$ (in dashed) traversed by a bi-infinite path in the dual graph (green), giving rise to an extremal state at low temperature.}
\label{Fig.lpqs}
\end{figure}

On hyperbolic lattices, Series and Sinai~\cite{series1990ising} exhibit an {uncountable family of mutually singular Gibbs states} indexed by the geodesics of $\HH_2$, see \Cref{Fig.lpqs}. The authors conjecture these states to be extremal. {The result {of Series and Sinai} holds for the ferromagnetic {\em n.n.}~Ising model defined on the vertices of the Cayley graph of any finitely generated co-compact group of isometries of $\HH_2$, see~\cite[Theorem 1]{series1990ising}.} 
More recently, Gandolfo--Ruiz--Shlosman~\cite{GRS15}, building on  Series--Sinai's paper, construct an uncountable family of inhomogeneous Gibbs states with multiple rigid interfaces called the \emph{Millefeuilles}. \\

In the present paper we adapt the techniques of \cite{ckln}  to {exhibit} an {uncountable} family of infinite-volume Gibbs measures, whose existence and extremality is proven by means of Peierls estimates and excess energy around local ground states. 

\noindent

\medskip
\noindent

\section{Results}

{
Let $\mathcal{L}_{p,q}$ be the tiling of the hyperbolic plane $\mathbb{H}_2$ such that each face is a {regular} $p$-gon and each vertex has degree $q$.
 
Our main result is a sufficient condition for well-definiteness and extremality of the low temperature Gibbs {measures} of the Ising model on $\mathcal{L}_{p,q}$ obtained via weak limits. These weak limits are obtained by imposing boundary conditions which have a sufficiently sparse set of frustrated bonds, but a priori {no symmetry of the underlining} $\mathcal{L}_{p,q}$ tilings{, whence {inhomogeneous}}. The key point is {an extension of the excess energy Lemma of~\cite{ckln} to {{non-amenable graphs} of bounded degree, which include in particular the graph given by the 1-skeleton of the hyperbolic lattices $\mathcal{L}_{p,q}$}, which we shall refer henceforth to simply by $\mathcal{L}_{p,q}$ by a slight abuse of notation. {Informally, this sufficient condition {measures} the excess energy {with respect to a reference configuration} in terms of the isoperimetric constant of the graph. 
\smallskip
	
For $H \subset G$ a finite subgraph of $G$, denoting by $\partial H$ the usual external boundary of $H$, the \dfn{isoperimetric constant} of $G$ is defined by 
$$
{\rm IC}_{G} \overset{\rm def}{=}\inf \left\{\frac{|\partial \gamma|}{|\gamma|}\sut \gamma \text{ a finite and non-empty subgraph of } G \right\}.
$$ 
Its explicit value on $\mathcal{L}_{p,q}$ has been famously determined by H\"aggstr\"om--Jonasson--Lyons~\cite[Theorem 4.1]{HJL2002}, and we will denote it by $\IC$ in this case.

The Ising model on a locally finite graph $G=(V,E)$ is defined as follows : take $\Lambda\subset V$ be a finite subset, and $\omega$ a fixed element in $\{-1,+1\}^V$. Define the Hamiltonian $H:\{-1,+1\}^\Lambda\to\R$ as
$$H^\omega_\Lambda(\sigma)=-\sum_{\{i,j\}\in E\cap\Lambda}\sigma_i\sigma_j-\sum_{\{i,j\}\in E\cap\partial\Lambda}\sigma_i\omega_j$$
At inverse temperature $\beta>0$, the Ising model in $\Lambda$ with boundary condition $\omega$ is the probability measure $\mu_\Lambda^\omega$  on $\Omega_\Lambda=\{-1,1\}^\Lambda$ proportional to $\exp(-\beta H^\omega_\Lambda(\sigma))$.
The set of infinite volume Gibbs measures of the Ising model at inverse temperature $\beta$ is the set of probability measures $\mu$  on $\Omega_\Lambda=\{-1,1\}^V$ such that the DLR equations hold:
$$\mathcal G_\beta=\{ \mu \in\mathcal M_1(\{-1,1\}^V) : \mu=\int d\mu(\omega)\mu^\omega_{\Lambda}\text{ for any finite  }\Lambda\subset V\}$$
The set $\mathcal G_\beta$ is a Choquet simplex~\cite{georgii}, and we denote its extremal elements by $\rm ex \mathcal G_\beta$. 
\newline
Let $\delta_{\rm max}\in\mathbb N$, and 
$$\Omega_{\rm GS}(\delta_{\rm max}):=\Big\{\omega\in\{-1,1\}^V : \forall i \in V, 
\sum_{\substack{j\in V\\i \sim j}}\mathbf{1}_{\omega_i\neq\omega_j}\leq \delta_{\rm max}\Big\}$$ 
be the set of (infinite) spin configurations such that there are at most $\delta_{\rm max}$ frustrated edges emanating from any vertex. We call these {{$\delta_{\rm max}$-inhomogeneous}} configurations, see \Cref{def.kd} for details in the case of $\mathcal{L}_{p,q}$. Then the following holds:
}
\begin{theorem}[{\textsc{{Sufficient condition for extremality}}}]\label{thm.speis}

Let ${\rm IC}_G$ be the isoperimetric constant of the {connected, transitive and locally finite graph $G$}. If the following sparsity condition holds
\begin{equation}\label{sparsity condition}
\delta_{\rm max} < \frac{1}{2}\; {\rm IC}_G
\end{equation}
then, at sufficiently low temperature, for any $\omega \in \Omega_{\rm GS}(\delta_{\rm max})$, the Gibbs measure $\mu^{\omega}$ of the Ising model on $G$, obtained as weak limit with boundary condition $\omega$, is {well-defined} and {extremal}.
\end{theorem}
\noindent
{Note that, for Theorem 2.1 to be useful, it is enough to get a (uniform and strictly positive) lower bound on ${\rm IC}_G$, which is however a non-trivial task generally.
When $G$ is  given by the 1-skeleton of $\mathcal{L}_{p,q}$, the isoperimetric constant $\IC$ is explicitely known~\cite{HJL2002} and the following corollary holds:
}
{
\begin{cor}[{\textsc{Uncountably many inhomogeneous extremal states {for the Ising model on $\mathcal L_{p,q}$}}}]\label{cor.umies}
\LC{}
\LC{For any $p>4$, $q\geq 3$ such that $\frac{1}{p}+\frac{1}{q}<\frac{1}{2}$ and ${\IC}>2$, there exists $\beta_0=\beta_{0}(p,q)$ such that, for any $\beta>\beta_{0}$, $\rm ex \mathcal G_\beta$ is {\bf uncountable}.}

\end{cor}
}

\begin{figure}
	\centering
	\begin{overpic}[width=.45\textwidth]{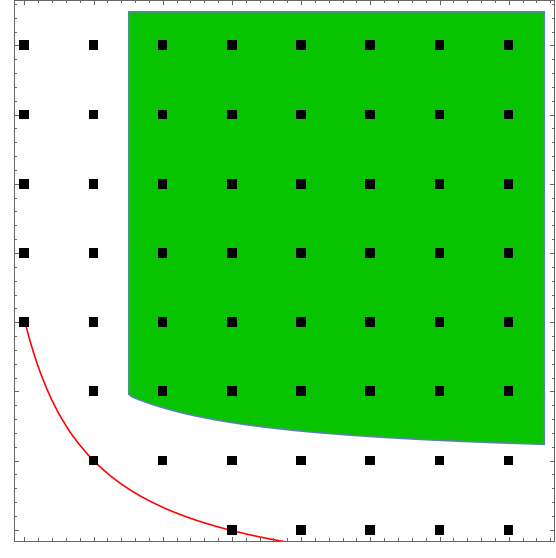}
	\put (4,-2) {\small $3$}
	\put (16,-2) {\small $4$}
	\put (28,-2) {\small $5$}
	\put (40.5,-2) {\small $6$}
	\put (52.7,-2) {\small $7$}
	\put (65.3,-2) {\small $8$}
	\put (77.5,-2) {\small $9$}
	\put (88.5,-2) {\small $10$}
	\put (-3,90.5) {\small $10$}
	\put (-2,78.5) {\small $9$}
	\put (-2,66) {\small $8$}
	\put (-2,54) {\small $7$}
	\put (-2,41.5) {\small $6$}
	\put (-2,29) {\small $5$}
	\put (-2,17) {\small $4$}
	\put (-2,4) {\small $3$}
	\put (102,2) {\large $p$}
	\put (1,102) {\large $q$}
	\end{overpic}
	\caption{Validity region (in green) of \Cref{cor.umies}. The red line is the curve $1/p+1/q=1/2$. The blue line is the curve $\IC>2$. 
	}
	\label{fig.validity-region}
\end{figure}

The corollary is a consequence of the observation that, in the validity region, uncountably many Dobrushin-like configurations ($+$ on one side of a separating line passing through the origin, $-$ on the other side) give rise to extremal states.  Indeed, interfaces of configurations on $\mathcal L_{p,q}$ can be represented as lines on the dual graph $\mathcal L_{q,p}$, see Figure \ref{Fig.lpqs}. Now, {following Moran~\cite{moran1997growth}}, let us draw the ``corona representation'' of the graph $\mathcal L_{q,p}$, see Figure \ref{fig.corona-representation}. The point is that for $q>4$, $\mathcal L_{q,p}$ strictly includes a union of $q$ trees of degree $(q-3)$ (which are glued at the origin). If one chooses an infinite branch in one of these trees (an end), and another one in one of the $q-3$ opposite trees (neighboring trees are forbidden), then we obtain a bi-infinite line $\gamma$ which has the required property: the sparsity condition \eqref{sparsity condition} is fulfilled for the Dobrushin configuration separated by $\gamma$. The quantity of such good interfaces $\gamma$ is uncountable by the usual Cantor set argument {for the ends of trees} {(see \emph{e.g.~}\cite[Ex. I.8 31]{BH})}.

\begin{figure}
	\centering
	\begin{overpic}[width=.45\textwidth]{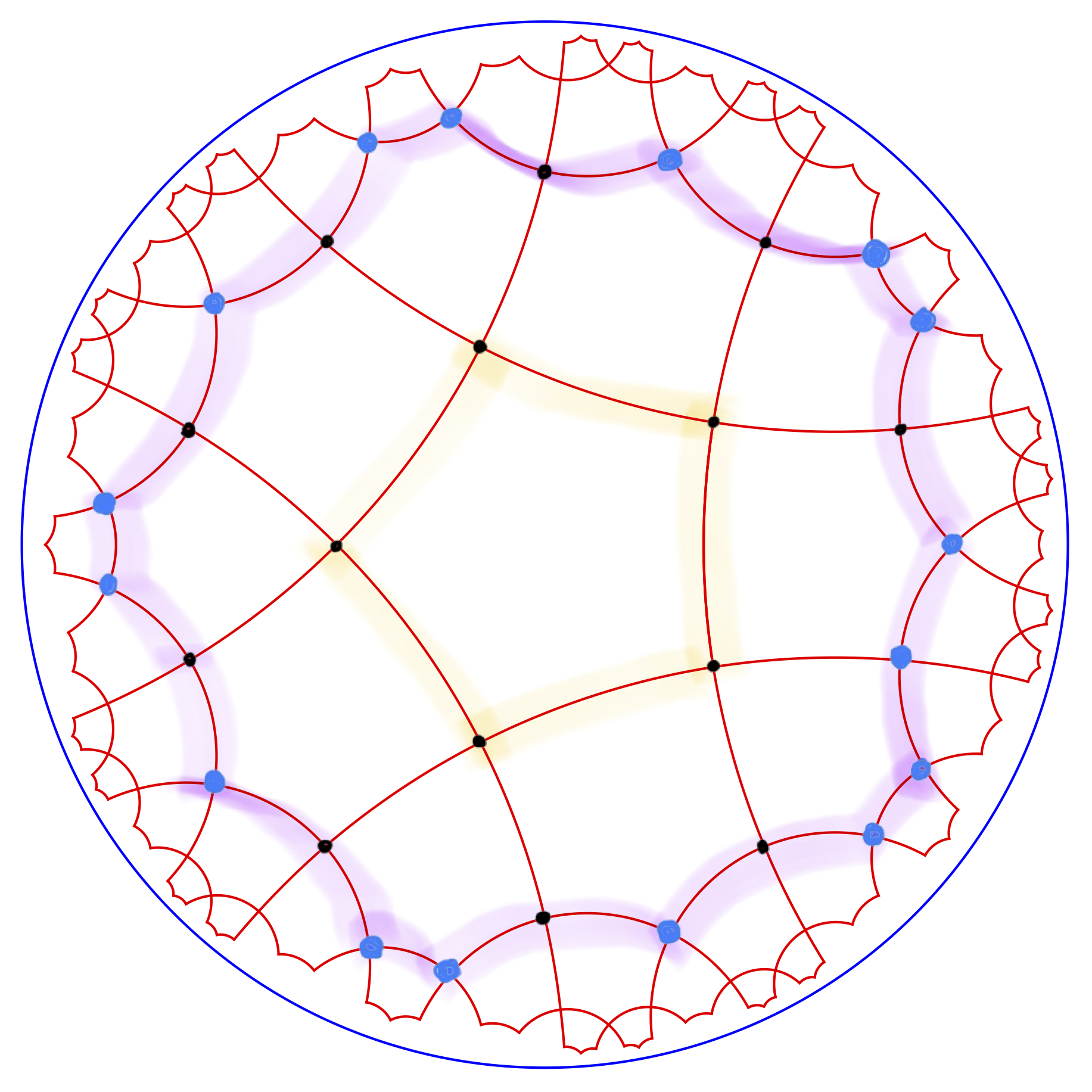}
	\put (52,50) {$\origin$}
\put (48.8,49.5) {{\tiny $\bullet$}}
	\end{overpic}
	\begin{overpic}[width=.45\textwidth]{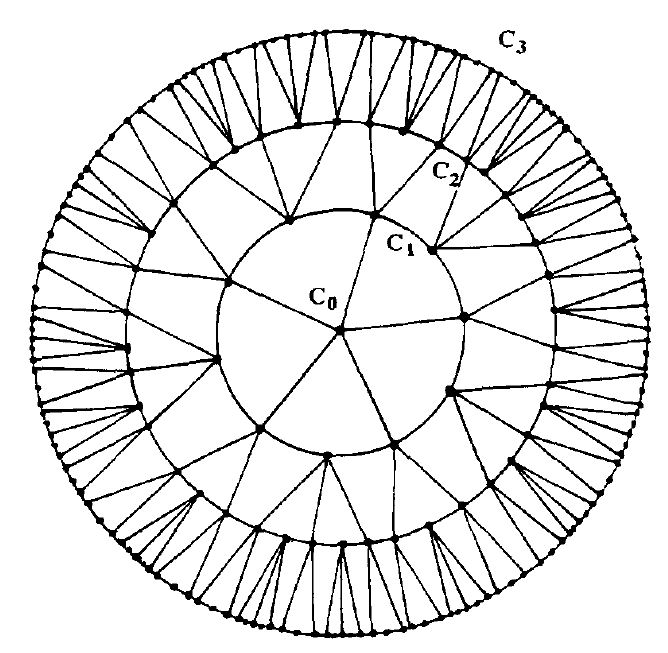}
	\end{overpic}
	\caption{LEFT : $\mathcal L_{4,5}$ in red, 
		RIGHT : Corona representation of $\mathcal L_{5,4}$,
		taken from \cite{moran1997growth}. }
	\label{fig.corona-representation}
\end{figure}

\medskip

Series and Sinai left open the question of extremality of their ``interface states''\footnote{{We quote here a statement in \cite[Page 896]{wuII}, which might suggest that the same authors proved the extremality of their states. However, the question is still open to the best of our knowledge.}} selected by a continuum geodesics $\gamma$ of $\HH_{2}$, see \cite[page 2, line -9]{series1990ising}.  {A natural question is to wonder if all the Series--Sinai states are extremal.} \Cref{thm.speis} above gives a sufficient condition for Series-Sinai states to be extremal: the {bi-infinite} geodesics $\gamma$ {of $\HH_2$}  should not cross more than $\frac\IC2 \leq \frac{q}{2}-1$ edges emanating from any vertex of $\mathcal{L}_{p,q}$. This condition can be nicely reformulated in terms of billiard (hyperbolic) trajectories in a regular hyperbolic $q$-gone, by folding $\gamma$ in a given fundamental domain via the reflection law. The sufficient condition is that the biliard trajectory corresponding to $\gamma$ in the Dirichlet domain of $\mathcal{L}_{p,q}$ should never bounce more than $\frac{q}{2}-1$ times onto two {adjacent} faces. 
A lemma due to Mostow \cite{Mostow} (see also \cite[Lemma 3.43]{K}) allows us to conjecture that uncountably many Series-Sinai states are extremal. Indeed, a true geodesic of $\mathbb H_2$ is close to each of the Dobrushin states constructed  by means of quasi-geodesics in Corollary \ref{cor.umies}. We leave the full proof as well as reformulations of this nice question in different domains of mathematics for a subsequent work.

%



\bigskip

\noindent \textbf{Structure of the paper}. \Cref{subsec.model} and \Cref{sec:def} contain the definition of the model and some key notions.
\Cref{sec.prsies}, which is a generalisation of the approach of~\cite{ckln}, is devoted to the proof of \Cref{thm.speis}. A key novelty {of this paper is} the generalisation of the ``excess energy lemma'' of \cite{ckln} to rather general graphs, see~\Cref{lem.ee}, in terms of their isoperimetric constant. Section \ref{sec.ld} focuses on Lobatchevsky planes $\mathcal L_{p,q}$, for which the isoperimetric constant is known, and allows to deduce the  uncountability of the set of extremal Gibbs states at low temperature.
\medskip

\section{Proofs}\label{sec.proofs}

\subsection{\newALN{Ising model}}\label{subsec.model}
\noindent
Let $G=(V,E)$ be a simple, locally finite connected graph.


\medskip

\noindent
{\bf Measurable structure}. We consider Ising spins: the single--spin space is $\Omega_{0}=\{-1,+1\}$, endowed with the $\sigma$-algebra given by the power set $\mathcal{P}(\{-1,+1\})$ and the a priori measure $\rho_{0}=\frac{1}{2}\left(\delta_{-1}+\delta_{+1} \right)$.

The (infinite-volume) configuration space is $\Omega=\Omega_{0}^{V}$, with events $\mathcal{F}=\left(\mathcal{P}(\{-1,+1\})\right)^{\otimes V}$ and a priori (infinite-volume, product) measure $\rho= \left(\rho_{0}\right)^{\otimes V}$. 

\medskip
\noindent
{\bf Microscopic finite-volume Hamiltonian}. For a finite subset of vertices $\Lambda \Subset V$, a configuration $\sigma \in \Omega$ and a boundary condition $\omega$, we consider the following \emph{ferromagnetic}, \emph{n.n.}~Ising Hamiltonian
$$
{
H^{G}_{\Lambda}(\sigma \mid \omega) \overset{\rm def}{=}  \sum_{\substack{ {v} \sim {w} \\ {v},{w} \in \Lambda}} \mathbf{1}_{\sigma_{{v}}\neq \sigma_{{w}}} + \sum_{\substack{ {v} \sim{u} \\ {v} \in \Lambda, \; {u} \in \Lambda^{c}}} \mathbf{1}_{\sigma_{{v}}\neq \omega_{{u}}} \; .}
$$




\medskip
\noindent
\subsection{Some definitions}\label{sec:def}
Let $G=(V,E)$ be a locally finite transitive graph.
\begin{defn}[\textsc{Contour wrt $\sigma^0$ and compatibility}]
Let $\sigma^0\in \Omega$. A contour is a  {connected} subgraph $\gamma \subset G$ s.t.~$\sigma_\gamma=1-\sigma^0_\gamma$. Two contours $\gamma,\gamma'$ are {compatible} if and only if $\gamma \bigtriangleup \gamma' = \gamma \cup \gamma'$, where $\bigtriangleup$ denotes the symmetric difference. We shall denote this compatibility relation by $\gamma \sim \gamma'$.
\end{defn}

For two configurations $\sigma, \eta \in \Omega$, we shall denote by $V(\sigma \bigtriangleup \eta) \subset V$ the set of sites at which $\sigma$ and $\eta$ disagree, which might be finite or infinite. Thus $\sigma \bigtriangleup \eta$ is a union of contours.

{
\noindent 
 Then we can provide the following definition:
\begin{defn}[\textsc{$\delta_{\rm max}$-inhomogeneous configurations}]\label{def.kd}
For a configuration $\sigma \in \Omega$, let $G_b(\sigma)$ be the subgraph of $G$ induced by the set of \emph{broken bonds} in $\sigma$, i.e.~$E_b(\sigma)=\{e=({v},{w}) \in E \sut\sigma_{v} \cdot \sigma_{w} =-1 \}$ (with the obvious vertex set $V_b(\sigma)$). \newline
Let 
$$
\delta_{\rm broken}(\sigma)=\max_{v \in V_b(\sigma)} {\rm deg}(v).
$$
 Let $\delta_{\rm max}\in \mathbb{N}$. Then
 $$
 \Omega_{\rm GS}(\delta_{\rm max})=\{\sigma^0 \in \Omega \sut \delta_{\rm broken}(\sigma^0) \leq \delta_{\rm max} \} )\subset \Omega 
 $$
 is called the set of $\delta_{\rm max}$-inhomogeneous configurations.
\end{defn}
}
\newpage









\subsection{Proof of Theorem 2.1}\label{sec.prsies}

Inspired by the work~\cite{ckln} about finite spin models on regular trees, we shall now prove the following:
\begin{lem}[\textsc{Excess energy lemma}]\label{lem.ee}
	Let $G$ be some connected graph. Let $\gamma$ be a finite contour wrt to a given configuration $\sigma^0 \in \Omega$ and let $\sigma=\sigma_\gamma \sigma^0_{\gamma^c}$ be the configuration coinciding with $\sigma$ in $\gamma$ and with $\sigma^0$ in $\gamma^{c}$. Then the following holds
	$$
	H(\sigma)-H(\sigma^0) \geq  \left({\rm IC}_G -2 \delta_{\rm broken}(\sigma_{0}) \right)|\gamma| \; ,
	$$
	\noindent
	where 
	${\rm IC}_G$
	is the isoperimetric constant of $G$.
\end{lem}

\rproof
First, by definition,
$$
H(\sigma)-H(\sigma^0)=\sum_{\substack{{v} \in \gamma \\ {w} \in \gamma^c \\ {v} \sim {w}}}\left(\mathbf{1}_{\sigma_{v} \neq \sigma_{w}}-\mathbf{1}_{\sigma^0_{v} \neq \sigma^0_{w}} \right) \; .
$$
Since the two configurations $\sigma$ and $\sigma^0$ differ precisely at {the contour} $\gamma$ and {by definition} ${v} \in \gamma$, ${w}\in \gamma^c$, we have $\sigma_{v}=\sigma_{w}$ if and only if $\sigma^0_{v}\neq \sigma_{w}^0$. Equivalently,  $\mathbf{1}_{\sigma_{v} \neq \sigma_{w}}=1-\mathbf{1}_{\sigma^0_{v} \neq \sigma^0_{w}}$. Thus 
$$
\sum_{\substack{{v} \in \gamma \\ {w} \in \gamma^c \\ {v} \sim {w}}}\left(\mathbf{1}_{\sigma_{v} \neq \sigma_{w}}-\mathbf{1}_{\sigma^0_{v} \neq \sigma^0_{w}} \right)=\sum_{\substack{{v} \in \gamma \\ {w} \in \gamma^c \\ {v} \sim {w}}}\left(1-  2 \cdot  \mathbf{1}_{\sigma^0_{v} \neq \sigma^0_{w}}\right) = |\partial \gamma| - 2 \sum_{\substack{{v} \in \gamma \\ {w} \in \gamma^c \\ {v} \sim w}}\mathbf{1}_{\sigma^0_{v} \neq \sigma^0_{w}} \; .
$$
By definition of ${\rm IC}_G$ we have a {lower bound} on the first term as $|\partial \gamma|\geq {\rm IC}_G |\gamma|$ and an {upper bound} on the second term by $2 \delta_{\rm broken}(\sigma_{0}) |\gamma|$ (as in~\cite[Proof of Lemma 2]{ckln}). This concludes the proof.
\Qed

\begin{remark}  In the case where $G=\mathcal L_{p,q}$, we have $\IC=(q-2)\sqrt{1-\frac{4}{(p-2)(q-2)}}$, see ~\cite[Theorem 4.1]{HJL2002}.
	Thus for any fixed $q\geq 3$, we have $\lim_{p\to \infty}\IC=(q-2)$ which recovers \cite[Lemma 2]{ckln} (to match the notations one should put $d=q-1$ and $u=U=1$ in \cite[Eq.~3.9]{ckln}). Also, $\IC=c_{p,q}\cdot \sqrt{\frac{1}{2}-\frac{1}{p}-\frac{1}{q}}$ (for some strictly positive constant $c_{p,q}$), which equals zero only for the three Euclidean lattices.
\end{remark}

Now let $\delta_{\rm max} \in \mathbb{N}$ be fixed. The proof of Theorem \ref{thm.speis}, i.e.~{well-definition and extremality} of $\mu^\omega$ for $\omega\in\Omega_{GS}(\delta_{\rm max})$ such that $\delta_{\rm max}<{\rm IC}_G/2$ follows now from Lemma \ref{lem.ee} via convergence of the cluster expansion and asymptotic decorrelation of polymer type events when $\beta$ is large enough, as in \cite{ckln}. For this we need the graph $G$ to be connected, transitive (to perform the weak limit), non-amenable (${\rm IC}_G>0$), and locally finite ($\delta_{\rm max}<\infty$). However, to actually have a non-empty sufficient condition for extremality, one needs to at least have a positive lower bound on the Cheeger constant, and this task can be highly non-trivial.
Note that in \cite{ckln} the argument {also} extends to a wide range of $SOS$-like models on trees, thanks to the fact that contours have no interior {on trees}. Here we stick to the Ising model in order to be able to perform the usual ``spin-flip {trick}'' when proving Peierls bounds. 
\subsection{Proof of Corollary 2.2}\label{sec.ld}
\subsubsection{Lobatchevsky planes}
 The lattice structure {we consider} is provided by a regular hyperbolic tiling of {the hyperbolic plane} {$(\mathbb{H}_{2},\mathrm{d}_{\mathbb{H}_{2}})$} in which tiles are congruent, regular geodesic polygons with $p$ sides of hyperbolic length equal to 1, meeting at vertices with degree $q$. Each such tiling is indexed by the Schläfli symbol $\{p,q\}$, for integers $p,q \geq 3$ s.t.~$\frac{1}{p}+\frac{1}{q}<\frac{1}{2}$ (see~\cite{grunbaum1987tilings}). With a slight abuse of notation we will  denote the family of infinite $q$-regular metric graphs associated to the $\{p,q\}$ tilings by $\mathcal{L}_{p,q}=(V_{p,q},E_{p,q},F_{p,q})$.

\noindent
Elements of the vertex set $V_{p,q}$ are denoted by {$v,w,\ldots$}. Two vertices {$v$, $w$} are said to be nearest neighbors, denoted $v \sim w$ (or abbreviated \emph{n.n.}) if and only if {$v$} and {$w$} are incident to the same edge $e \in E_{p,q}$. {We denote by ${\rm d}_{gr}$ the graph distance on $\mathcal{L}_{p,q}$.} Equivalently, due to our embedding convention, {$v$} and {$w$} are \emph{n.n.}~if and only if $\mathrm{d}_{\mathbb{H}_{2}}(v,w)=1$. 

{ Let $G$ be a subgraph of $(V_{p,q},E_{p,q})$. We write $|G|$ for the number of vertices in $G$ ($\equiv|V(G)|$). The set of \emph{emanating edges} from $G$, which is a subset of $E_{p,q}$, is defined by  $\partial G=\{e=(v,w) \in E_{p,q} \sut v \in V(G) \; \text{and} \;  w \not\in V(G) \}$. The \emph{inner boundary} of $G$, which is a subset of $V_{p,q}$, is defined by $\partial_{in}G=\{v \in V(G) \sut (v,w) \in E_{p,q} \; \text{and} \;   w \not \in V(G)  \}$.} {For $G_1,G_2$ two given subgraphs of $\mathcal{L}_{p,q}$, their distance is defined by ${\rm d}_{gr}(G_1,G_2)=\min_{\substack{v \in G_1\\ w \in G_2}}{\rm d}_{gr}(v,w)$.}

{For $r\geq 0$ and a vertex $v$, the combinatorial ball of radius $r$ centered at $v$ is defined as follows $B(r,v)\overset{\rm def}{=}\{w \in V_{p,q} \sut {\rm d}_{gr}(v,w)\leq r\}$. Analogously, the combinatorial sphere of radius $r$ centered at $v$ is defined by $\partial B(r,v)\overset{\rm def}{=}\{w \in V_{p,q} \sut {\rm d}_{gr}(v,w)= r\}$.}

\subsubsection{Corona representation}

{$\mathcal{L}_{p,q}$ admits a \emph{layer decomposition} (see~\cite[Appendix A]{mertensmoore2017}) which allows to get exact formulas for the number of vertices, faces and edges at or within a given layer from any given point $z \in \mathbb{H}_{2}$ via a transfer matrix method. This method has been discovered independently by Rietman--Nienhuis--Oitmaa~\cite[Eq.~2.14]{RNO} and Moran~\cite[Pag. 159]{moran1997growth}, who has also generalized it to {other} homogeneous tilings of $\mathbb{H}_2$. It has been extensively used by the first author together with Vanessa Jacquier and Wioletta M.~Ruszel to exhibit the set of finite shape with minimal perimeter, whose exhaustion of $\mathcal{L}_{p,q}$ actually realize $\IC$, see~\cite{DJR}. One consequence {of this layer decomposition} is {an exact formula for} the number of vertices at or within a given layer. This number grows exponentially in the layer number with a growth rate given by the largest eigenvalue of the transfer matrix.

{We will recall here {these} results for face--centered lattices $\mathcal{L}_{p,q}$ in which $z$ is the barycenter of the fundamental polygon and it is put at the origin $\origin$ of $\mathbb{H}_{2}$ (see Figure~\ref{Fig.lpqs}). Also, we will recall the layer decomposition for vertices of $\mathcal{L}_{p,q}$; the layer decomposition for faces and edges of $\mathcal{L}_{p,q}$ is analogous but not needed for our current purposes.}
}

Layer zero $L_{0}$ is the set of vertices $v \in V_{p,q}$ of the face containing $\origin$. The first layer $L_1$ is the set of vertices of the faces sharing a vertex or edge with the central face in $L_{0}$, which are not in $L_{0}$; and so on \emph{ad libitum} {(see \Cref{fig.ld})}. Let $S_{n;p,q}=\# \lbrace  v \in V_{p,q} \sut v \in L_n  \rbrace $ 
and $B_{n;p,q}=\sum_{i=0}^{n}S_{i;p,q}$ be respectively the cardinal of layer $n$ and the total number of vertices within layer $n$ (layer included). 

\begin{figure}
\centering
\begin{overpic}[width=.45\textwidth]{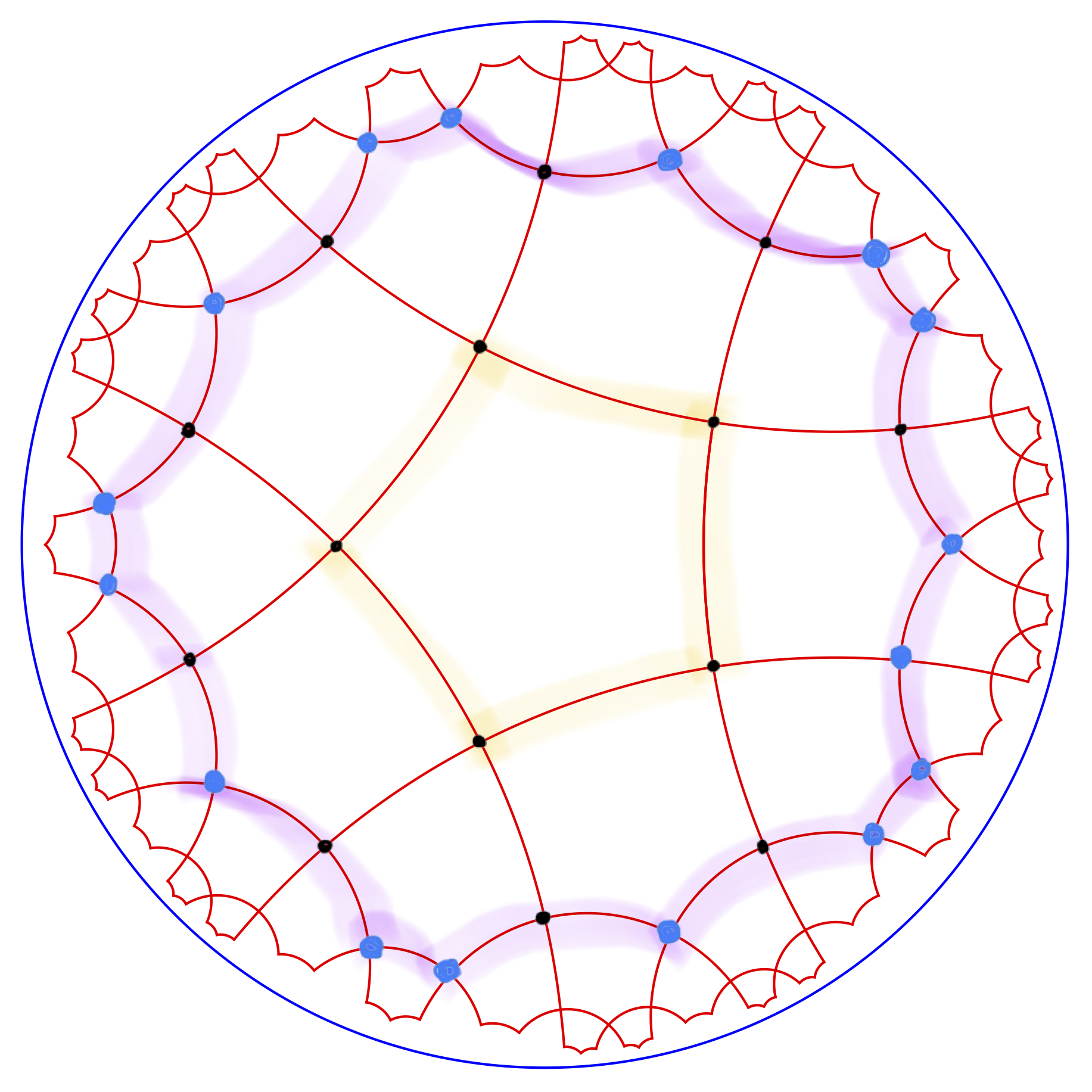}
 \put (55,66) {$L_{0}$}
 \put (18.5,78) { $L_{1}$}
 \put (82,3) { $\mathcal{L}_{5,4}$}
  \put (47.3,48.7) { $\bullet$}
   \put (52,52) { $\origin$}
\end{overpic}
\caption{{Layer decomposition for $\mathcal{L}_{5,4}$ {(first two layers)}. In each layer, vertices are of two kind (in blue and black here) allowing counting by induction, see~\cite{RNO,moran1997growth} for details.}}
\label{fig.ld}
\end{figure}

\begin{figure}
	\centering
	\begin{overpic}[width=.45\textwidth]{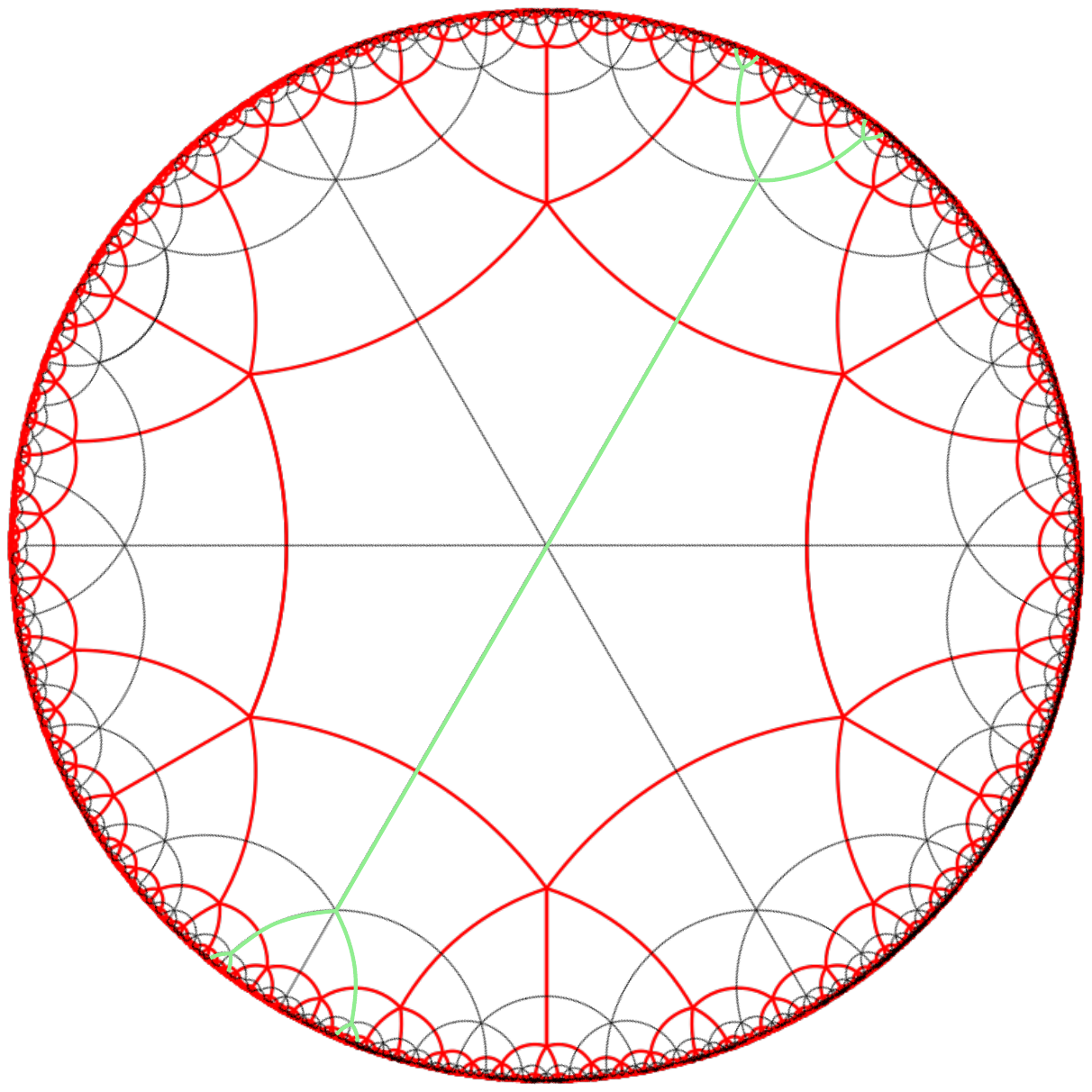}
	\put (35.5,52) {$\origin$}
\put (36,49.5) {{\tiny $\bullet$}}
	\end{overpic}
	\begin{overpic}[width=0.45\textwidth]{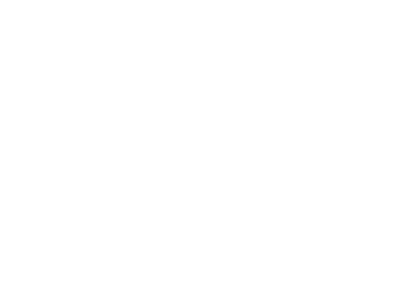}
     \put(10,50){\includegraphics[width=0.45\textwidth]{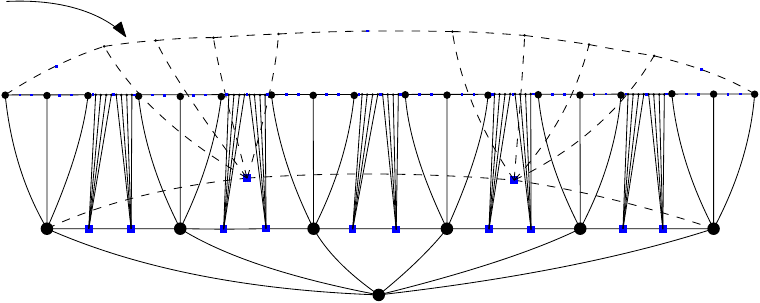}}  
     \put (60,45) {$\origin$}
  \end{overpic}
	\caption{LEFT: $\mathcal L_{5,6}$ in red, its dual $\mathcal L_{6,5}$ in black, and two of the embedded ternary trees {(glued at $\origin$)} in green. RIGHT : The corona representation of $\mathcal L_{6,5}$. }
	\label{fig.ld2}
\end{figure}

The corona representation of $\mathcal{L}_{p,q}$ goes as follows, see Figure \ref{Fig.lpqs}:

\begin{itemize}
	\item From the origin, draw $p$ outgoing edges, and their $p$ end-vertices in black (first generation).
	\item Draw a circle around them (first circle, consisting of $p$ vertices and $p-1$ edges).
	\item Add $q-3$ blue vertices on each edge of the circle (triangle$\rightarrow q$-gone\footnote{To complete the $q$-gone from the 3 existing edges.}).
	\item Draw $p-3$ outgoing edges from the black vertices of the circle, and $p-2$ outgoing edges from the blue vertices of the circle, and draw all their end-vertices in black (second generation).
	\item Draw a new circle around this second generation of vertices (second circle).
	\item Add $q-3$ blue vertices on each edge whose endpoints are attached to the same vertex of the first circle (triangle$\rightarrow q$-gone), and $q-4$ blue vertices on each edge whose endpoints are attached to different vertices (4-gone$\rightarrow q$-gone\footnote{To complete the $q$-gone with the existing 4 edges.}) of the first circle.
	\item Draw $p-3$ outgoing edges from the black vertices of the second circle, 
						$p-2$ outgoing edges from the blue vertices of the second circle, 
						and all their end-vertices in black (third generation).
	\item Repeat the procedure \emph{ad libitum}. 
\end{itemize}

Note that for $p>4$, the corona representation draws a union of $p$ disjoint trees of degree $p-3$, which are glued at the origin, and follow the edges linking one circle to the next one. 

\subsubsection{Building interfaces ensuring the sparsity condition}
Interfaces of the configurations on the graph $\mathcal L_{p,q}$ can be represented as lines on the graph $\mathcal L_{q,p}$, see Figure \ref{Fig.lpqs}. We will {now} build uncountably many Dobrushin configurations which fulfill the sparsity condition \eqref{sparsity condition}. In the corona representation $\mathcal L_{q,p}$, we just proved that for $q>4$, there is a union of $q$ trees of degree $q-3$ (which are glued at the origin). 

Choose an infinite branch in one of these trees, and another one in one of the $q-3$ opposite trees (neighboring trees are forbidden), then we obtain a bi-infinite line $\gamma$ which do not cross more than one edge emanating from each vertex of the primal graph $\mathcal L_{p,q}$. Let thus $\delta_{\rm max}$ be equal to 1. Put $+$ spins on one side of $\gamma$ and $-$ spins on the other side. This forms a Dobrushin configuration $\omega^\pm_\gamma$, which belongs to $\Omega_{GS}^{p,q}(1)$. 

Thus, whenever $\IC>2$, the sparsity condition \eqref{sparsity condition} holds, and by Theorem \ref{thm.speis}  the infinite-volume Gibbs measure, obtained as weak limit with boundary condition $\omega^\pm_\gamma$, is {well-defined} and {extremal}. The quantity of such good interfaces $\gamma$ is uncountable by the usual Cantor set argument {(see~\emph{e.g.} \cite[Example 8.11.5]{BH})}. This finishes the proof of Corollary \ref{cor.umies}.\qed

{ 
}

\vfill

\noindent \textbf{Acknowledgements}.
The authors thank Fran\c cois Dahmani, Federica Fanoni, Greg McShane and Ana Rechtman for enlightening discussions. {A.L.N.~thanks Senya Shlosman for having introduced him to these beautiful problems}. M.D'A.~thanks Alessio Lerose for pointing out reference~\cite{BPR2020} to him. 
M.D'A.~is grateful to the Institut Fourier in Grenoble for excellent working conditions  of two invitations (January and June 2024). 
The work of M.A.~is supported by the ERC Consolidator Grant SuperGRandMA (Grant No.~101087572).

\bibliographystyle{acm}

\clearpage
\bibliography{biblio}

\end{document}